\documentclass[12pt]{article}
\usepackage{amsmath}
\usepackage{latexsym}
\usepackage{amsfonts}
\def\Bbb{\mathbb}

\def\eea{\end{eqnarray*}}

\newtheorem{main}{Theorem}

\newtheorem{defn}{Definition}
\newtheorem{thm}{Theorem}
\newtheorem{prop}[thm]{Proposition}

\newtheorem{lem}[thm]{Lemma}

\newenvironment{proof}{\medskip \noindent
{\bf Proof.}}{\hfill \rule{.5em}{1em}
\\}

\def\RR{{\mathbb R}}

\begin{document}

\title{Perelman's  Invariant, Ricci Flow, and  the Yamabe 
 Invariants of Smooth Manifolds}

\author{Kazuo Akutagawa, Masashi Ishida,   and Claude LeBrun\thanks{Supported 
in part by  NSF grant DMS-0604735.}  
  }

\date{October 2, 2006\\
Revised: October 17, 2006}
\maketitle

\begin{abstract}
In his study of Ricci flow, Perelman introduced a smooth-manifold
invariant  called $\bar{\lambda}$. We show here  that, 
for completely elementary reasons,  this invariant  simply equals the 
Yamabe invariant, alias  the sigma constant, whenever the latter is non-positive. 
On the other hand,  the Perelman invariant  just  equals $+\infty$  whenever
the Yamabe invariant is positive.  
\end{abstract}

Let $M$ be a smooth compact manifold of dimension $n\geq 3$. 
Perelman's celebrated work on Ricci flow \cite{perelman1,perelman2} 
led him to consider the
functional which associates to every Riemannian metric $g$
the least eigenvalue $\lambda_g$ of the elliptic operator $4 \Delta_g+s_g$,
where $s_g$ denotes the scalar curvature of $g$, and 
$\Delta = d^*d= - \nabla\cdot\nabla $ is the positive-spectrum  Laplace-Beltrami operator 
associated with  $g$. In other words, $\lambda_g$ can be expressed in terms
of Raleigh quotients as 
$$
\lambda_g = \inf_{u} \frac{\int_M \left[ s_gu^2 + 4 |\nabla u|^2 \right] d\mu}{\int_M u^2d\mu}
$$
where  the infimum is  taken over all smooth, 
 real-valued functions $u$ on $M$. 
 
 One of Perelman's remarkable observations is that the scale-invariant quantity
 $\lambda_gV_g^{2/n}$ is non-decreasing under  the  Ricci flow, where 
 $V_g= \int_Md\mu_g$ denotes the total volume of $(M,g)$. This led him to 
 consider the  differential-topological invariant obtained by taking the 
 supremum of this quantity over the space of all Riemannian metrics \cite{perelman2,lott}:

\begin{defn}
Let $M$ be a smooth compact $n$-manifold, $n\geq 3$. Perelman's 
$\bar{\lambda}$ invariant of $M$ is defined to be 
$$
\bar{\lambda}(M)= \sup_g \lambda_gV_g^{2/n}
$$
where 
the supremum is taken over all smooth metrics $g$ on $M$. 
\end{defn}

The r\^ole of the scalar curvature and Laplacian  in defining $\lambda_g$
might immediately make one wonder whether this invariant might somehow
be related to the Yamabe problem. Recall that, as was conjectured
by Yamabe \cite{yam}, and later proved by Trudinger, Aubin, and Schoen  
\cite{aubyam,aubis,lp,rick,trud}, 
every conformal class on any smooth compact manifold contains a metric of 
constant scalar curvature. If 
$M$ is a smooth compact manifold of dimension $n\geq 3$, and if 
$$\gamma=[g]=\{ vg ~|~v: M\to {\Bbb R}^+\},$$ 
is the conformal class of an arbitrary metric, such a metric $\hat{g}$ 
can in fact be constructed by minimizing the normalized total-scalar-curvature 
functional  
$$
\hat{g}\mapsto  \frac{\int_M 
s_{\hat{g}}~d\mu_{\hat{g}}}{\left(\int_M 
d\mu_{\hat{g}}\right)^{\frac{n-2}{n}}},
$$
among all metrics conformal to $g$. Indeed, by setting
$\hat{g} = u^{4/(n-2)}g$, this expression can be rewritten as
$$
\frac{\int_M 
s_{\hat{g}}~d\mu_{\hat{g}}}{\left(\int_M 
d\mu_{\hat{g}}\right)^{\frac{n-2}{n}}}= 
\frac{\int_M\left[ s_gu^2 +
4 \frac{n-1}{n-2}|\nabla u|^2\right] d\mu_g}{\left(\int_M  u^{2n/(n-2)}d\mu_g\right)^{(n-2)/n}},
$$
and the proof proceeds by showing that there is a smooth positive function 
$u$ which minimizes the right-hand  expression. In particular,
 each conformal class $\gamma$ has an  associated number 
$Y_{\gamma}$, called its 
{\em Yamabe constant}, obtained by setting 
$$Y_{\gamma} = \inf_{g\in \gamma}  \frac{\int_M 
s_{{g}}~d\mu_{{g}}}{\left(\int_M 
d\mu_{{g}}\right)^{\frac{n-2}{n}}}~;$$
 the content of the Trudinger-Aubin-Schoen theorem is exactly that 
this number is actually realized as the constant scalar curvature of some
unit-volume metric in $\gamma$. A constant-scalar-curvature metric 
of this type is called a {\em Yamabe minimizer}. It is not
difficult to show that any Riemannian metric with $s=\mbox{const} \leq 0$
is  a Yamabe minimizer, and that, moreover, the   Yamabe
minimizer $g\in \gamma$ 
is unique (up to  constant rescaling) 
whenever  $Y_\gamma \leq 0$. The situation is much more complicated
when $Y_\gamma > 0$, but  it is still not difficult  to see that if $g$
is a metric for which $s$ has a fixed sign (positive, zero, or negative) everywhere on $M$, 
then 
this sign necessarily agrees with that of the number $Y_{[g]}$.  

Yamabe's work was  apparently
 motivated by the hope of constructing Einstein metrics via 
 a variational approach. This idea eventually led Kobayashi \cite{okob} and Schoen \cite{sch}
 to independently introduce the smooth manifold invariant 
$${\mathcal Y}(M) = \sup_{\gamma} Y_{\gamma} = 
\sup_{\gamma}\inf_{{g}\in \gamma } \frac{\int_M
s_{g}~d\mu_{g}}{\left(\int_M 
d\mu_{g}\right)^{\frac{n-2}{n}}}.$$
By construction, this is a diffeomorphism 
invariant of $M$, and  is now commonly known as the 
  {\em Yamabe invariant} of $M$; note, however, that 
  Schoen called ${\mathcal Y}(M)$  the {\em sigma constant},
 and that this terminology  is still preferred by some authors. 
  Notice that ${\mathcal Y}(M) \leq 0$ iff 
   $M$ does not admit
metrics of positive scalar curvature, and that, when this happens, 
${\mathcal Y}(M)$ 
is simply the   supremum of 
the scalar curvatures of unit-volume 
constant-scalar-curvature metrics on $M$.

The fact that there is some fundamental relation between 
the Yamabe invariant ${\mathcal Y}(M)$ and Perelman's $\bar{\lambda}$ invariant
was probably first  pointed out by Anderson \cite{ander34}. More
recently, an e-print by Fang and Zhang \cite{fang}  computed the Perelman invariant
for a large class of $4$-manifolds where the Yamabe invariants had already been 
computed by the present authors \cite{lno,il2} and others \cite{jp2,jp3}, and,
as was later emphasized by Kotschick  \cite{dkot}, 
their answers  exactly agree with those previously discovered   in the Yamabe
case. The point of this brief note is to observe that this was no mere matter
of coincidence:

\begin{main}
Suppose that $M$ is a smooth compact $n$-manifold, $n\geq 3$. Then 
$$\bar{\lambda}(M) = \begin{cases}
     {\mathcal Y}(M) & \text{ if  } {\mathcal Y}(M) \leq 0, \\
     +\infty  & \text{ if  } {\mathcal Y}(M) >  0.
\end{cases}
$$
\end{main}

In fact,  this will follow easily once we clearly  understand the  behavior of 
$\lambda V^{2/n}$ on  each individual conformal class.

\begin{prop} \label{yupyup} 
Suppose that $\gamma$ is a conformal class on $M$ which does not contain
a metric of positive scalar curvature. Then 
$$Y_\gamma = \sup_{g\in \gamma} \lambda_g V^{2/n}_g$$
\end{prop}
\begin{proof}
Let $g\in \gamma$, and let $\hat{g}= u^{4/(n-2)}g$  be the Yamabe minimizer
in $\gamma$. Then 
$$0\geq Y_\gamma = \frac{\int_M\left[ s_gu^2 +
4 \frac{n-1}{n-2}|\nabla u|^2\right] d\mu_g}{\left(\int_M  u^{2n/(n-2)}d\mu_g\right)^{(n-2)/n}}.$$
Thus 
\begin{eqnarray*}
\lambda_g \int u^2 d\mu &\leq& \int \left[ su^2 + 4|\nabla u|^2 \right]d\mu \\
&\leq&  \int \left[ su^2 + 4\frac{n-1}{n-2} |\nabla u|^2 \right]d\mu \\
&=& Y_\gamma \left(\int  u^{2n/(n-2)}d\mu\right)^{(n-2)/n} \\
&\leq& Y_\gamma V^{-2/n} \int u^2 d\mu 
\end{eqnarray*}
where, since $Y_\gamma\leq 0$,
the last step is an the application of the H\"older inequality
$$\int f_1f_2 ~d\mu \leq \left(\int |f_1|^pd\mu \right)^{1/p} \left(\int |f_2|^qd\mu \right)^{1/q}, 
~~~\frac{1}{p}+ \frac{1}{q}=1,$$
with $f_1=1$, $f_2=u^2$, $p= n/2$, and $q=n/(n-2)$. 
Moreover, equality holds precisely when $u$ is constant  --- which is to say, 
precisely when 
$g$ has constant scalar curvature.

Since this shows that 
$$ \lambda_g  V^{2/n} \leq Y_\gamma$$
for every $g\in \gamma$, and  since equality occurs if $g$ is the
Yamabe minimizer, it follows   that 
$$
\sup_{g\in \gamma} \lambda_g V^{2/n}_g = Y_\gamma , 
$$
exactly as claimed. 
\end{proof}

We now need  make only one more simple observation:

\begin{lem} \label{aha} 
If $M$ carries a metric with $s > 0$, then $\bar{\lambda}(M) = +\infty$. 
\end{lem}
\begin{proof}
Given such a manifold $M$ and any smooth non-constant function
$f: M\to \RR$,  Kobayashi \cite{okob} has shown that 
there exists a unit-volume  metric  on $M$ with 
$s=f$. In particular, given any real number $L$, there is 
 a unit-volume metric $g_L$ on $M$ with $s > L$ everywhere. 
 But for such a metric, $\lambda > L$ and $V=1$. Thus, taking $L\to \infty$, 
 $\bar{\lambda}(M) = \sup_g \lambda_gV_g^{2/n}= + \infty$. 
\end{proof}

Theorem A is now follows immediately. Indeed, if ${\mathcal Y}(M) >0$,
$M$ admits a metric with $s > 0$, and Lemma \ref{aha} therefore tells
us that $\bar{\lambda}(M) = +\infty$. Otherwise, no conformal class
contains a metric of positive scalar curvature, and Proposition 
\ref{yupyup} therefore tells us that each constant-scalar-curvature metric
maximizes $\lambda V^{2/n}$ in its conformal class. Given an arbitrary 
maximizing sequence  $\hat{g}_j$ for  $\lambda V^{2/n}$,  we may
thus, by conformal rescaling, 
construct  a new maximizing sequence  ${g}_j$ consisting of 
unit-volume constant-scalar-curvature metrics. But for any such sequence, 
 the numbers $s_{g_j}$ may be viewed either as  $\{ \lambda_{g_j}V_{g_j}^{2/n}\}$ or as  $\{Y_{[g_j]}\}$. Thus the  suprema over the space of all Riemannian metrics of $Y_{[g]}$ and 
$\lambda_g V_g^{2/n}$ must 
precisely coincide. 

Now, there is a substantial literature \cite{il2,lno,lric,jp2,jp3,petyun} concerning manifolds of
 non-positive Yamabe invariant, and the exact value of the
invariant is moreover known for large numbers of such manifolds. By virtue of 
Theorem A, all of these facts   about ${\mathcal Y}(M)$  
 may therefore immediately be interpreted 
as  instead pertaining to    $\bar{\lambda}(M)$. 

On the other hand, we have also seen that the Perelman invariant 
jumps to positive infinity whenever there is a metric of 
positive scalar curvature. By contrast, the Yamabe invariant 
always remains finite; indeed, one of Aubin's fundamental contributions
to the theory of the Yamabe problem is the fact that 
${\mathcal Y}(M) \leq {\mathcal Y}(S^n)$ for any smooth compact 
$n$-manifold $M$. This systematic discrepancy fundamentally reflects 
the fact that constant-scalar-curvature metrics are generally not  
minimizers in the positive case \cite{okob,rick}, and 
that, moreover, constant-scalar-curvature metrics of arbitrarily 
high energy exist in profusion  \cite{pollack} in this setting. Thus,
while we have found  a fascinating  link between 
Perelman's invariant and the Yamabe problem, 
this observation actually teaches us nothing at all  precisely in the
regime where the Yamabe problem is the most technically difficult and,
in some respects, still the most poorly understood.

\vfill

{\footnotesize 
\noindent
{Kazuo Akutagawa, 
{Dept.\  Mathematics, Tokyo Univ.\   of Science, Noda  278-8510,  
Japan}\\{\sc e-mail}: akutagawa\underline{~}kazuo@ma.noda.tus.ac.jp}\\
{Masashi Ishida,
{Department of Mathematics,  
Sophia University, 
 Tokyo  102-8554, Japan}\\
{\sc e-mail}: ishida@mm.sophia.ac.jp}
 \\ {
Claude LeBrun,
Department of Mathematics, SUNY, 
Stony Brook, NY 11794-3651, USA\hfill 
{\sc e-mail}: claude@math.sunysb.edu}}

\end{document}